\documentclass[pdflatex,sn-vancouver-num]{sn-jnl}

\usepackage{graphicx}%
\usepackage{multirow}%
\usepackage{amsmath,amssymb,amsfonts}%
\usepackage{amsthm}%
\usepackage{mathrsfs}%
\usepackage[title]{appendix}%
\usepackage{xcolor}%
\usepackage{textcomp}%
\usepackage{manyfoot}%
\usepackage{booktabs}%
\usepackage{algorithm}%
\usepackage{algorithmicx}%
\usepackage{algpseudocode}%
\usepackage{listings}%
\usepackage{bm}
\usepackage{bbm}

\theoremstyle{thmstyleone}%
\newtheorem{theorem}{Theorem}
\newtheorem{proposition}[theorem]{Proposition}%

\theoremstyle{thmstyletwo}%

\theoremstyle{thmstylethree}%

\def\ignore#1{}
\def\cH{\mathcal{H}}
\def\cX{\mathcal{X}}
\def\cA{\mathcal{A}}
\def\tcA{\widetilde{\mathcal{A}}}
\def\law{\mathcal{L}}
\def\cF{\mathcal{F}}
\def\cG{\mathcal{G}}
\def\cZ{\mathcal{Z}}
\def\bbP{\mathbb{P}}
\def\bbX{\mathbb{X}}
\def\bbG{\mathbb{G}}
\def\bP{\mathbf{P}}
\def\bQ{\mathbf{Q}}
\def\Def{\ :=\ }
\def\TV{{\mathrm{TV}}}
\def\BW{{\mathrm{BW}}}
\def\re{\mathbb{R}}
\def\ex{\mathbb{E}}

\def\Eq{\ =\ }
\def\Le{\ \le\ }

\def\Po{\mathrm{Po}}
\def\Be{\mathrm{Be}}
\def\la{\lambda}
\def\bi{\mathbf{i}}
\def\bk{\mathbf{k}}
\def\cI{\mathcal{I}}
\def\hI{\widehat{I}}
\def\nin{\noindent}
\def\bbA{\mathbb{A}}
\def\Pop{\mathbf{PP}}
\def\bla{\boldsymbol{\la}}

\begin{document}

\title[Stein's method in network analysis]{Stein's method in network analysis}

\author[1]{\fnm{A. D.} \sur{Barbour}}\email{a.d.barbour@math.uzh.ch}

\author[2]{\fnm{Adrian} \sur{Fischer}}\email{adrian.fischer@stats.ox.ac.uk}

\author[2]{\fnm{Gesine} \sur{Reinert}}\email{reinert@stats.ox.ac.uk}

\affil[1]{\orgdiv{Institut f\"ur Mathematik}, \orgname{Universit\"at Z\"urich}, 
\orgaddress{\street{Winterthurerstrasse 190}, \city{8057 Z\"urich},   \country{Switzerland}}}

\affil[2]{\orgdiv{Department of Statistics}, \orgname{University of Oxford}, \orgaddress{\street{24-29 St Giles'}, \city{Oxford}, \postcode{OX1 3LB}, \state{Oxfordshire}, \country{United Kingdom}}}


\abstract{ The paper consists of a brief survey of the use of Stein's method in network analysis.
Topics covered include normal and Poisson approximation of graph--based statistics, approximating an
exponential random graph by a Bernoulli model, and comparison of different random geometric graph models.}




\maketitle

\section{Introduction}
\ignore{When analysing networks, summary statistics, for example based on the distribution of degrees or on local clustering, are often used. Understanding the theoretical behaviour of these summary statistics under different models can inform model choice and assessing the significance of deviations from what may be expected. However, even in simple models such of that of a Bernoulli random graph, with independent and identically distributed edges, the distribution of the number of triangles already does not have a closed form. Hence suitable approximations are required. Moreover, many social networks are often fairly small, and hence the approximation error may not be negligible. To address these issues, Stein's method has been established as a tool of choice.}

To determine whether a network might plausibly have arisen from a particular network model, it is usual to compare values
of a number of summary statistics with what might be expected according to the model. Typical statistics include those based on
counts of small subgraphs, on the empirical degree distribution, and on local measures of clustering. Yet, even for the
simplest models, the exact distributions of such statistics are not known.  It may then be useful to establish approximations
to their distributions, with some measure of their accuracy.  Stein's method is a good tool for this.

Stein's method for a probability distribution~$\bP$ on a space~$\cX$ starts with a \textit{Stein operator}~$\mathcal{A}$
for~$\bP$, such that, if $Y \sim \bP$, then
\begin{align} \label{stein_identity}
    \mathbb{E}[\mathcal{A}f(Y)] \Eq 0 \qquad \mbox{for all } f \in \cF,
\end{align}
where~$\cF$ is a suitably large class of functions $f\colon \cX \to \re$.  
If the  \textit{Stein equation}
\begin{align} \label{stein_equation}
    h(x) - \mathbb{E}[h(Y)] \Eq \mathcal{A}f(x)
\end{align}
has a solution $f := f_h \in \cF$ for each $h \in \cH$,  where~$\cH$ is a determining class for~$\bP$, then, if~$X$ is any random element of~$\cX$ such that $\ex[\cA f(X)]$ is uniformly small for all $f \in \cF$, it follows that $\ex\{h(X) - \mathbb{E}[h(Y)]\}$ is uniformly small for all~$h \in \cH$. This means that the distribution $\law(X)$
of~$X$ is close to~$\bP$, in the sense that $d_{\cH}(\law(X),\bP)$ is small, where the distance~$d_{\cH}$ is defined by
$$ 
   d_\cH(\law(U),\law(V)) \Def \sup_{h \in \mathcal{H}} \vert \ex h(U) - \ex h(V)  \vert. 
$$
For instance, if $\cX = \mathbb{R}^d$, 
choosing
\[
 \cH \Eq \mathcal{H}_{\TV} \Def \{h\colon\mathbb{R}^d \rightarrow \mathbb{R}; \ 
 h(t) = \mathbf{1}\{t \in B\} ,
 B \mbox{ a Borel set in } \mathbb{R}^d \}
\]
gives $D_{\cH} = d_{\TV}$, the total variation distance;  and choosing
\[
 \cH \Eq \mathcal{H}_{\BW} \Def \{h\colon\mathbb{R}^d \rightarrow \mathbb{R};\ |h(t_1)| \le 1,\ 
 \vert h(t_1) - h(t_2) \vert \leq \Vert t_1-t_2 \Vert \text{ for all } t_1,t_2 \in \mathbb{R}^d \}
\]
gives $D_{\cH} = d_{\BW}$, the bounded Wasserstein distance.

There are two standard approaches for
obtaining a Stein operator~$\cA$.  The first is the \textit{generator approach} \cite{barbour1988stein,barbour1990stein,gotze1991rate}, in which~$\mathcal{A}$ is the generator of an ergodic
Markov process~$Z_{\bP}$, whose stationary distribution is~$\bP$.  If the Markov process is tractable, it
directly yields a solution to the Stein equation~\eqref{stein_equation}.
The second is the \textit{density approach} \cite{stein1986approximate,ley2013stein,mijoule2023stein}, in which, 
for $\cX = \re$,
$\mathcal{A}f=(pf)'/p$, where~$p$ is the density of the distribution~$\bbP$. 
For $\cX = \re^d$, with $d \ge 2$, there are a number of options:  see \cite{mijoule2023stein}.
If solutions to~\eqref{stein_equation} are known, it remains only to bound~$\ex[\cA f(X)]$ for $f \in \cF$,
something typically achieved using Taylor expansion and coupling techniques, as illustrated in
\cite{barbour1992poisson,chen2010normal,ross2011fundamentals}
and \cite{anastasiou2023stein}.
Stein's method has been used in a wide variety of different settings, resulting in a vast amount of literature\footnote{\url{https://sites.google.com/site/steinsmethod/articles}}. 
Here, we focus on random graphs.

\ignore{Just a note: the density approach here is phrased only for one-dimensional distributions. For higher dimension there are several options, such as taking the gradient operator $\nabla (pf)/p$; \cite{mijoule2023stein} has more  details.}

We denote by $G=(V,A)$ a graph on~$n = |V|$ vertices with adjacency matrix $A = (A_{u,v})$, where $A_{u,u} = 0$ and,
for $u \ne v$,
$A_{u,v}=A_{v,u}=1$ if there is an edge between $u$ and $v$, and $A_{u,v}=A_{v,u}=0$ otherwise. In a random graph, $A$ is  random. In the Bernoulli random graph model $\mathcal{G}(n,p)$, the edge indicators are independent,
with $\bbP(A_{u,v}=A_{v,u}=1) =p$ for all $u \ne v \in V$. In an Erd\H{o}s--R\'{e}nyi mixture graph, each vertex belongs to exactly one class,
and the independent edge indicators have probabilities depending only on the {\it classes\/} of their end points.  
We discuss the approximation of count statistics in these Bernoulli random graph models in Section~\ref{section_bernoulli}.
Exponential random graph models are popular families having dependent
edge indicators, and their approximation by Bernoulli models is illustrated in Section~\ref{section_ergm}. The spatial random graph models in Section~\ref{section_spatial_graphs} have conditionally independent edge probabilities, given the positions of the vertices in a underlying space.

\section{Bernoulli models} \label{section_bernoulli}
Motif counts are summary statistics frequently used in network analysis. They can be represented in the form $X = \sum_{\bi \in \cI}I_\bi$, where~$\bi$ is a (possibly ordered) set of vertices,
$\cI$ is the collection of all such~$\bi$, and~$I_\bi$ is a function of the entries of~$A$ with indices in~$\bi$.
Standard examples include the numbers of edges, triangles, $k$-stars and isolated vertices. For instance, for $G = (V,A)$,
\begin{gather*}
    \text{\#triangles } \Eq \sum_{ \{u,v,w\} \subset V} a_{u,v} a_{v,w} a_{w,u}, \qquad 
    \text{\#2-stars } \Eq \frac{1}{2} \sum_{u,v,w \in V, v \neq w} a_{u,v} a_{u,w},\\
    \text{\#isolated vertices } \Eq \sum_{v \in V} \prod_{u \in V} (1-a_{u,v}).
\end{gather*}
When~$A$ is the adjacency matrix of
a random graph with independent edge indicators (thus following what we call here a {\it Bernoulli model\/}), the distributions of such statistics (suitably normalised) can often be well approximated by Poisson or normal distributions. 
Stein's method can be used to assess the accuracy of such approximations; see, for example, \cite[Chapters 7 and~9]{barbour2026networks}.
\ignore{
converge to a certain distribution 
as the number of vertices tends to infinity. With Stein's method as introduced in the introduction one can determine conditions on the the random graph model as to when this convergence to a certain target distribution holds and also quantify the approximation by an explicit bound on an integral probability metric. Common approximating distribution in this setting are the Poisson and normal distribution.
}

\subsection{Poisson approximation}\label{Poisson_approximation}
When a count~$X$ is typically not very large, a
Poisson distribution $\Po(\lambda)$ with expectation $\lambda = \ex X$ may offer a plausible approximation. 
A suitable Stein operator for~$\Po(\la)$ is given by \citep{chen1975poisson}
\begin{align}\label{Poisson_Stein_operator}
    \mathcal{A}f(x) \Eq \cA_\la f(x) \Eq \lambda f(x+1)-xf(x), \qquad x= 0,1,\ldots\,,
\end{align}
and can be obtained from the density approach by using a discrete derivative $\Delta f(x) = f(x+1) - f(x)$. 
Alternatively, writing $f(x) = g(x) - g(x-1)$, the Stein operator~$\cA_\la$ can be recognized as 
equivalent to the 
generator~$\tcA_\la$ of an immigration-death process~$Z_\la$ with immigration rate~$\lambda $ and unit 
per capita death rate:
\begin{align}\label{Poisson_Stein_operator2}
    \tcA_\la g(x) \Eq \lambda ( g(x+1) - g(x) )+ x ( g(x -1)- g(x) ), \qquad x= 0,1,\ldots\,.
\end{align} 
The solution~$f_h$ to \eqref{stein_equation} for 
$h \in \mathcal{H}_{\TV}$, with~$\cA$ as in~\eqref{Poisson_Stein_operator}, satisfies
\begin{align*}
    \Vert f_h \Vert_\infty \Le \min\{1 , \lambda^{-1/2} \}\quad\mbox{and} \quad \Vert \Delta f_h \Vert_\infty 
       \Le \min\{1 , \lambda^{-1} \} \quad 
    \text{for all }  h \in \mathcal{H}_{\TV},
\end{align*}
where we write $\Vert \cdot \Vert_\infty$ for the supremum norm. First, consider $X = \sum_{i=1}^m I_i$, 
where $I_i \sim \Be(p_i)$ are independent Bernoulli random variables with expectation~$p_i$. 
Then, for $\lambda = \mathbb{E}[X] = \sum_{i=1}^m p_i$, using elementary manipulations 
\citep[(1.12)--(1.15)]{barbour1992poisson},
we have
$$
     \mathbb{E}[\mathcal{A}_\la f(X)]  
     \Eq \sum_{i=1}^m p_i^2 \mathbb{E}[\Delta f(X-I_i +1)].
$$ 
Applying the bound on~$\Delta f_h$ from above, it follows that
\begin{align}\label{indep-bnd}
    d_{\TV}(\law(X),\Po(\lambda)) \Le \min\{1, \lambda^{-1} \} \sum_{i=1}^m p_i^{2}.
\end{align}
As an example, consider the number of edges in an Erd\H os--R\'enyi mixture graph $\mathcal{G}_K(\mathbf{n},P)$ with~$K$ types,
where $\mathbf{n}=\{n_1, \ldots, n_K \}$ are the numbers of vertices of each type and $P \in (0,1)^{K \times K}$ is the 
symmetric matrix of edge probabilities. Then the number of {\it edges\/} can indeed be represented as a sum~$X$ of independent random
variables, as above. The bound~\eqref{indep-bnd} gives
\begin{align} \label{bernoulli_mixture_bound}
      d_{\TV}(\law(X),\Po(\lambda)) \Le \min\{1, \lambda^{-1} \} \bigg( \sum_{k=1}^K \binom{n_k}{2} p_{k,k}^2 + \sum_{k=1}^{K-1} \sum_{l=k+1}^K n_k n_lp_{k,l}^{2} \bigg),
\end{align}
where 
$$
     \lambda = \sum_{k=1}^K \binom{n_k}{2} p_{k,k} + \sum_{k=1}^{K-1} \sum_{l=k+1}^K n_k n_lp_{k,l}\,.
$$ 
The right-hand side of \eqref{bernoulli_mixture_bound} is bounded by $\max_{1 \leq k \leq l \leq K} p_{k,l}$,
which is frequently small. 

For counts $X = \sum_{\bi \in \cI}I_\bi$ of motifs containing more than one edge, the indicators are no longer mutually independent, because $I_\bi$ and~$I_{\bi'}$ may not be independent if $\bi$ and~$\bi'$ overlap, and the bound~\eqref{indep-bnd} cannot be used.
However, if the occurrence of a particular~$I_\bi$ has little effect on the total~$X$, the dependence may be weak enough for
Poisson approximation still to be plausible.  More precisely, suppose that, for each $\bi \in \cI$, there is a
decomposition
\ignore{
A natural question is if the approach extends to dependent random variables $I_i$. For this purpose we consider cases where we have some sort of local dependence, i.e.\ $I_i$ is independent from most but not all the other random variables involved in the sum. We therefore assume that
}
$X = X_\bi + Z_\bi + I_\bi$, where $X_\bi$ and~$I_\bi$ are independent.
\ignore{
If~$Z_\bi$ is not too large compared to~$X$.
For instance, if~$X$ counts triangles, with~$i$ representing a particular triangle,  $X_i$ might consist of all
triangles having no edges overlapping~$i$, in which case 
contains the random variables that depend on $I_i$. 
}
Then it follows from \cite[(2.2)]{chen1975poisson} (see also \cite[Theorem~1.A]{barbour1992poisson}) that, for $\la = \ex X  = \sum_{\bi\in\cI} p_\bi$,
\begin{align} \label{poisson_bernoulli_local_bound}
    d_{\TV}(\law(X),\Po(\lambda)) \Le \min\{1, \lambda^{-1} \} 
        \sum_{\bi\in\cI} p_\bi^{2}+ p_\bi \mathbb{E}[Z_\bi]+ \mathbb{E}[I_\bi Z_\bi].
\end{align}
For the number of triangles in the Bernoulli random graph $(V,A) \sim \mathcal{G}(n,p)$, the collection~$\cI$ is the set of all
triples $\bi = \{u,v,w\}$ of distinct vertices in~$V$, and $I_\bi = A_{u,v} A_{v,w} A_{w,u}$.
\ignore{
\begin{align} \label{sum_triangles_bernoulli}
     X = \sum_{ \{u,v,w\} \subset V} a_{u,v} a_{v,w} a_{w,u}.
\end{align}
}
Then $\lambda = \binom{n}{3}p^3$, and $I_\bi$ is independent of the collection 
$\bbX_\bi := \{I_{\bi'}\colon \bi' \in \cI_\bi\}$, where~$\cI_\bi$ consists of all triples~$\bi'$
such that $\bi'$ has at most one vertex in common with $\bi$, because the edges defining~$\bbX_\bi$ are disjoint from
those joining vertices of~$\bi$. 
\ignore{
Hence pick $3$ vertices $\star = \{u^{\star},v^{\star},w^{\star}\}$ and write
\begin{align*}
    X = X_{\star} + Z_{\star} + I_{\star},
\end{align*}
where $X_{\star}$ sums over all $\{u,v,w\}$ that have at most one vertex in common with $\star$ and $Z_{\star}$ 
sums up over the rest (except $I_{\star}$ itself).
}
Hence, writing $X_\bi = \sum_{\bi' \in \cI_\bi}$ and $Z_\bi = X - X_\bi - I_\bi$, it follows that
\begin{align*}
    \mathbb{E}[Z_{\bi}] \leq 3(n-3)p^3, \qquad \mathbb{E}[I_{\bi} Z_{\bi}] \leq 3(n-3)p^5,
\end{align*}
and the inequality \eqref{poisson_bernoulli_local_bound} gives
\begin{align*}
    d_{\TV}(\law(X),\Po(\lambda)) \Le \min\{1, \lambda \} (3n-8)p^2(1+p).
\end{align*}
The latter bound is small if $np^2 \ll 1$.

There is also a way to handle global dependence between random variables~$I_\bi$ when the dependence is weak, 
in that the value of a particular~$I_\bi$ has only a small effect on the the sum $X$. 
If $X_\bi^{\star}$ is constructed on the same probability space as~$X$, and
$\law(X_\bi^{\star}) = \law\bigl( \sum_{\bi', \bi' \neq \bi} I_{\bi'}\,|\,I_i = 1\bigr)$, then one can show
\citep[Theorem~1.B]{barbour1992poisson} that
\begin{align} \label{poisson_coupling}
     d_{\TV}(\law(X),\Po(\lambda)) \Le \min\{1, \lambda^{-1} \} \sum_{\bi} p_\bi \mathbb{E}\vert X- X_\bi^{\star} \vert.
\end{align}
As an example, consider the number of isolated vertices in $(V,A) \sim \mathcal{G}(n,p)$. Here, the random variables
$I_v = \prod_{u \in V} (1-A_{u,v})$ all depend on each other, but a given vertex being isolated only has a small 
effect on the probabilities that other vertices are isolated.  Indeed, defining
\begin{align*}
    X_v^{\star} \Def \sum_{u \in V, u \neq v} \prod_{w \in V, w \neq u,v} (1-A_{w,u})
\end{align*}
to be the number of vertices that are isolated in the restriction $G^v$ of~$G$ to $V\setminus\{v\}$,
we note that $X_v^{\star} \ge X - I_v$, {\it with equality if $I_v=1$},  and that 
\[
    \ex[X_v^{\star} - (X - I_v)] \Eq p\ex[X_v^{\star}] \Eq (n-1)p(1-p)^{n-2}.
\]
\ignore{
  $\ex\{X_v^{\star} - (X - I_v)\} = p\ex X_v^{\star} \le pe^{-(n-2)p}$ is small if~$p$
We construct a coupling $X_u^{\star}$ for a fixed vertex 
$u \in V$ as described above as follows by considering the graph without the vertex we conditioned on. Then,
\begin{align*}
    X_u^{\star} = \sum_{v \in V, v \neq u} \prod_{w \in V, w \neq u} (1-A_{w,v}).
\end{align*}
}
Hence, with $\lambda = \mathbb{E}[X] = n(1-p)^{n-1}$, it follows from \eqref{poisson_coupling} that
\begin{align*}
    d_{\TV}(\law(X),\Po(\lambda)) \Le (1-p)^{n-1} + (n-1)p(1-p)^{n-2} \Le (1+np)e^{-(n-2)p},
\end{align*}
small if $np \gg 1$.

\subsection{Normal approximation}
If $\bP = N(0,1)$, the density approach leads to the Stein operator
\begin{align}\label{Stein_normal_operator}
    \mathcal{A}f(x) \Eq xf(x)-f'(x).
\end{align}
Alternatively, $N(0,1)$ is the stationary distribution of the standard Ornstein--Uhlen\-beck process, and
the generator approach yields an equivalent Stein operator $\widetilde{\cA}g(x) := xg'(x) - g''(x)$.
For sufficiently regular~$h$, the solution to the Stein equation~\eqref{stein_equation},
with~$\cA$ as in~\eqref{Stein_normal_operator}, is given by
\begin{align*}
 f_h(x) &\Eq -e^{x^2/2} \int_{x}^{\infty} e^{-t^2/2} (h(t) - \mathbb{E}[h(Y)] )\, dt \\
 &\Eq e^{x^2/2} \int_{-\infty}^{x} e^{-t^2/2} (h(t) - \mathbb{E}[h(Y)] )\, dt
\end{align*}
\citep[(2.5)]{chen2010normal},
where $Y \sim N(0,1)$. 
\ignore{
Here, we consider the \textit{bounded Wasserstein distance} which is also an integral probability metric and 
obtained by the function class $\mathcal{H}_{BW}= \{h:\mathbb{R}^d \rightarrow \mathbb{R} \, \vert \, h 
\text{ lipschitz and } \Vert h \Vert_{\infty}, \Vert h' \Vert_{\infty} \leq 1 \} $, where $\Vert \cdot \Vert_{\infty} $ denotes the supremum norm. Similar to the Poisson case, 
}
For $h \in \mathcal{H}_{\BW}$, the class of functions defining the bounded Wasserstein distance,
the functions~$f_h$ satisfy
\begin{align}\label{normal-bnds}
    \Vert f_h \Vert_{\infty} \leq \sqrt{2 \pi}, \qquad \Vert f_h' \Vert_{\infty} \leq 4, 
         \qquad \Vert f_h'' \Vert_{\infty} \leq 4+ \sqrt{2 \pi}.
\end{align}

The following local approach is reminiscient of that leading to the Poisson local bound~\eqref{poisson_bernoulli_local_bound}.
We begin with $X := \sum_{i=1}^m I_i$, where the~$I_i$ do not necessarily have to have Bernoulli distributions;
write $\sigma^2 := \mathbb{V}[X]$ for the variance of $X$, and set $\hI_i := I_i - \ex[I_i]$.
Suppose that there is a {\it double decomposition\/} \cite{barbour1989central} of the following form.
First, for each $i=1, \ldots,m$, we can write $X = X_i + Z_i$, where~$X_i$ is independent of $I_i$.
Then, for each~$i$, we have $Z_i = \sum_{k=1}^{m_i} Z_{i,k}$, where, for each $1 \le k \le m_i$,
we can write $X_i = X_{i,k} + V_{i,k}$, with $X_{i,k}$ independent of $(I_i, Z_{i,k})$. 
\begin{theorem} \label{theorem_local_normal}
With the double decomposition above,  if all the $I_i$, $Z_{i,k}$ and~$V_{i,k}$ have finite third moments, then
    \begin{align*}
        d_{\BW}\bigl(&\law(\sigma^{-1}  (X- \mathbb{E}[X])),N(0,1)\bigr) \\
        &\Le  7\sigma^{-3} \, \sum_{i=1}^m \bigg( \tfrac12 \mathbb{E}[\vert \hI_i \vert Z_i^2] 
        + \sum_{k=1}^{m_i} \mathbb{E}[ \vert \hI_i Z_{i,k} V_{i,k} \vert ]
            + \mathbb{E}[ \vert \hI_i Z_{i,k}  \vert ]
        \mathbb{E}[ \vert Z_{i} +  V_{i,k} \vert ] \bigg).
    \end{align*}
\end{theorem}
\nin The proof of Theorem~\ref{theorem_local_normal} uses a Taylor expansion applied to $\ex[\cA f(X)]$
and the estimates~\eqref{normal-bnds}; see \cite[Theorem~9.3.1]{barbour2026networks}.

The number of triangles $X = \sum_{\bi\in\cI}I_\bi$ in a Bernoulli random graph $(V,A) \sim \mathcal{G}(n,p)$,
where $\cI$ is the set of all triples of distinct vertices and~$I_\bi$ is the indicator of the triangle on~$\bi$ 
being present, can also be approximated by a standard normal distribution. 
Here, as before, $m = |\cI| = \binom{n}{3}$, and~$I_\bi$ is independent of $X_\bi := \sum_{\bi'\in\cI_\bi}I_{\bi'}$,
where~$\cI_\bi$ consists of all triples such that~$\bi'$ has at most one vertex in common with $\bi$; thus we can take 
$Z_\bi = X - X_\bi$, and $Z_\bi = \sum_{\bk \in \cI \setminus\cI_\bi}Z_{\bi,\bk}$, with $Z_{\bi,\bk} = I_\bk$.
Then define $X_{\bi,\bk} := \sum_{\bi' \in \cI_\bi \cap \cI_\bk}I_{\bi'}$, independent of the pair $(I_\bi,I_\bk)$,
and set $V_{\bi,\bk} := X_{\bi} - X_{\bi,\bk}$.  This provides a suitable double decomposition, for use with
Theorem~\ref{theorem_local_normal}.  The variance $\sigma^2$ of~$X$ is given by
\ignore{
if they share at most one vertex. Therefore, we let $X_i$ to be the sum over all $I_j$ which are independent of $I_i$
and let $X_{i,k}$ to be the sum over all $I_j$ which are independent of $I_i$ and $I_k$ (hence $Z_{i,k} = I_k$).
furthermore, let $V_{i,k} = X_i - X_{i,k}$. Then $X_i$, $X_{i,k}$, $Z_{i,k}$ and $V_{i,k}$ satisfy the assumptions of Theorem \ref{theorem_local_normal}. Using the above decomposition of $X$ one can calculate the variance
}
\begin{align*}
    \sigma^2 \Eq \binom{n}{3} p^3\{(1-p^3) + 3(n-3)p^2 (1-p)\}.
\end{align*}
Bounds on  the expectations in the bound of Theorem~\ref{theorem_local_normal} can also be determined, in terms of $p$
and~$n$. It then follows from Theorem~\ref{theorem_local_normal} that there are constants $C_1$, $C_2$ and $C_3$,
independent of~$n$, such that
\begin{align*}
    d_{\BW}\bigl(\law(\sigma^{-1}  (X- \mathbb{E}[X])),N(0,1)\bigr) \Le \begin{cases}
         C_1 (np)^{-3/2}, &\ 0 < p \leq n^{-1/2}; \\
         C_2 (n\sqrt{p})^{-1}, &\ n^{-1/2} < p \leq  1/2; \\
         C_3 (n\sqrt{1-p})^{-1}, &\ 1/2 < p \leq 1.   \end{cases}
\end{align*}
The details are to be found in \cite[Example~7.3.3]{barbour2026networks}.  Note that, when $n^{-1} \ll p \ll n^{-1/2}$, both normal and Poisson approximations are reasonable. 

Another approach, known as \textit{size bias coupling}, can be seen as a generalisation of the coupling approach used above
for Poisson approximation. A random variable~$X^{(s)}$ is said to have the size bias distribution of a non-negative random variable~$X$ if, for any bounded real function~$f$,
\begin{align*}
    \mathbb{E}[f(X^{(s)})] \mathbb{E}[X] \Eq \mathbb{E}[Xf(X)].
\end{align*}
For example, if~$X\in \{ 0, 1, 2, \ldots\}$ then
$\mathbb{P}(X^{(s)} = j) = j \mathbb{P}(X = j) / \mathbb{E}[X] $. Proposition~\ref{theorem_size_bias_coupling}
shows how to size bias a sum of non-negative random variables $X = \sum_{i=1}^m Y_i$ with finite means. The proposition is implicit in the discussion following \cite[Lemma~2.1]{goldsteinrinott}; see also
\cite{Baldietal}.

\begin{proposition} \label{theorem_size_bias_coupling}
   Let $X = \sum_{j=1}^m Y_j$ be a sum of non-negative random variables, with $0 < \mathbb{E} [X] < \infty$,
   and let $F_j(\cdot;y) = \law\bigl(\sum_{i, i \neq j} Y_i \,|\, Y_j = y\bigr)$. Let $J$ be a random index with values in $1, \ldots, m$, with $\mathbb{P}(J = j) = \mathbb{E}[Y_j] / \mathbb{E}[X]$.
   Given $J=j$, let~$Y_j^{(s)}$ be sampled from the size biased distribution of~$Y_j$; then, if $Y_j^{(s)} = y$,
    sample~$X_j^{\star}$ from $F_j(\cdot;y)$. Then $X^{(s)} = X_J^{\star} + Y_J^{(s)}$ has the size bias distribution of~$X$.
\end{proposition}

\noindent In many applications, it is possible to use the recipe
of Proposition~\ref{theorem_size_bias_coupling} to construct~$X^{(s)}$ on the same probability
space as~$X$, in such a way that $X^{(s)}$ and~$X$ are close.
\ignore{
As their versions are slightly different, here is the proof of
Proposition~\ref{theorem_size_bias_coupling}.
For any  real function~$f$, by size biasing,
$\mathbb{E}[ f(X_j^{\star}(I_j^{(s)})+ I_j^{(s)} ] E[I_j] = E[I_j f(X_j^{\star}(I_j)+ I_j)]$,
and so
\begin{align*}
    \mathbb{E} f(X^{(s)}) \Eq& \sum_j \frac{ \mathbb{E}[I_j]}{\mathbb{E}[X]}\, \mathbb{E} [ f(X_j^{\star}(I_j^{(s)})+ I_j^{(s)} ]\\
    \Eq& \sum_j \frac{ \mathbb{E}[I_j]}{\mathbb{E}[X]}
       \frac{1}{\mathbb{E}[I_j]} \mathbb{E} \Bigl[ I_j f\Bigl(\sum_{i: i \ne j} I_i + I_j\Bigr) \Bigr] \Eq \mathbb{E} [X  f(X)]
\end{align*} as required.
}%
This motivates the following theorem;  see
\cite[Theorem~9.4.2]{barbour2026networks}.
\begin{theorem} \label{theorem_size_bias_bound}
 Let~$X$ be a non-negative random variable with finite variance~$\sigma^2$, and let~$X^{(s)}$ be constructed on the
 same probability space as~$X$, having the size bias distribution of~$X$;
 write $\nu^2 := \mathbb{E}[(X^{(s)}-X)^2]$. Let~$\mathscr{F}$ be any sub-sigma field with respect to the probability
 space on which $X$ and $X^{(s)}$ are defined.
 Then,
    \begin{align*}
        d_{\BW}\bigl(&\law(\sigma^{-1}  (X- \mathbb{E}[X])),N(0,1)\bigr) \\
        &\Le \sigma^{-2} \mathbb{E}[X] \bigg( 4 (\mathbb{V}[\mathbb{E}[X^{(s)}-X \, \vert \, \mathscr{F} ] ]  )^{1/2} +
        \frac{7\nu}{\sigma} \bigg( \frac{\nu}{2} + (\mathbb{E}[\mathbb{V}[X \, \vert \, \mathscr{F}]]^{1/2} \bigg) \bigg).
    \end{align*}
\end{theorem}
As an illustration, the number of vertices~$X$ of a given degree~$k$ in $(V,A) \sim \mathcal{G}(n,p)$ can be expressed as sum of Bernoulli
random variables,
\begin{align*}
    X \Eq \sum_{v \in V} \mathbf{1} \{\text{vertex } v \text{ has degree } k\}.
\end{align*}
Using the size bias distribution of~$X$ obtained from Proposition~\ref{theorem_size_bias_coupling}, it follows from
Theorem~\ref{theorem_size_bias_bound}  that $d_{\BW}\bigl(\law(\sigma^{-1}  (X- \mathbb{E}[X])),N(0,1)\bigr)
\le Cn^{-1/2}$, for a suitable constant~$C$, if $p = c/n$ for fixed~$c$;  details are to be found in
\cite[Example~9.4.4]{barbour2026networks}.
\ignore{
As the sub-sigma field $\mathscr{F}$ one can choose the sigma-field generated by the random variables $\{D_v \, \vert \, v \in V \setminus J \}$ where $D_v$ denotes the degree of vertex $v$ and $J$ is a uniform random variables in $V$ as introduced in the statement of Theorem \ref{theorem_size_bias_coupling}.}

We note that multivariate Poisson and normal approximation can also be carried out using Stein's method;
see \cite[Section~8.4 and Chapter~10]{barbour2026networks}.

\subsection{Compound Poisson approximation}

Counts of motifs other than triangles are used  to characterise networks. In \cite{topirceanu2016uncovering}, they are used to differentiate between online social networks like Facebook, Twitter, and Google--Plus. This study employs motifs of sizes 3 and~4, and counts them in collections of 50 ego-networks, of sizes 150 to 5000, from each of these social media platforms. They then record the percentage of times that each of the possible undirected motifs occur in each network group. For example, the triangle--whisker motif makes up 32.44\% of the motifs of size~4 in Facebook; for Twitter, the proportion is 27.33\%, and for Google-Plus it is 28.86\%. Some of the findings on motifs of size~4 are   summarised in Table \ref{tab:motifs}.

\begin{table}[]
    \centering
\begin{tabular}{c|c| c | c }
    platform & star & triangle--whisker motif & chain  \\ \hline
    Facebook  & 17.49 (low) & 32.44 (high) & 31.75 (high) \\
    Twitter & 22.50 (average-high) & 27.33 (average-high) & 30.43 (high)  \\
    Google-Plus & 31.48 (high) & 28.86 (high) &  21.34 (low) \\ \hline
\end{tabular}
\caption{Relative motif counts in ego-networks of three social media platforms; data and interpretation from \cite{topirceanu2016uncovering}}
    \label{tab:motifs}
\end{table}

To assess which differences are significant,  understanding the approximate count distributions under different
network models is useful.
Such distributions need not be approximately Poisson or normal. As an illustration, consider the triangle--whisker motif
in the Bernoulli random graph $\cG(n,p)$; for details, see \cite[Example 8.2.5]{barbour2026networks}.
If~$p$ is small, the distribution of the number of triangles may be close to Poisson,
but  each triangle may give rise to up to $3(n-3)$ triangle-whisker graphs; given that a triangle has occurred, each
possible extension to a triangle--whisker graph needs just one additional edge,
so that the expected number of such extensions is about~$3np$, which may be substantial.
Thus, one could think of occurrences appearing in clumps: each triangle can create a clump of up to $3(n-3)$
triangle--whisker graphs.

Such counts can often be well approximated by a compound Poisson distribution. Stein's method for compound Poisson
distributions is summarised for example in \cite{barbour2001compound}; results specialised to Bernoulli random graphs
can be found for example in \cite{stark2001compound}, and \cite{coulson2018compound} covers subgraph counts
in stochastic blockmodels.

\subsection{Assessing model fit}\label{modelfit}
Statistical conclusions drawn from a model rely on the model being appropriate for the data. To this purpose, \cite{fatima2025pure} developed tests for assessing the fit of a specified inhomogeneous random graph model in high dimensions. Inhomogeneous random graph models have independent edge indicators \citep{janson2007inhomogeneous}, without any structure 
being imposed on the edge probabilities.  
The tests are based on the \textit{Stein discrepancy} between probability distributions $\bP$ and $\bQ$, defined by
\begin{align}\label{Stein_discrepancy}
    S(\bP,\bQ) = \sup_{f \in \mathcal{H}} \vert \mathbb{E}[\mathcal{A}f(W)] \vert,
\end{align}
where $\mathcal{A}$ is a Stein operator with respect to $\bP$ and $W \sim \bQ$;  \cite{liuetal} contains further background and applications of Stein discrepancies. Taking the set $\mathcal{H}$ to be the unit ball in a reproducing kernel Hilbert space yields a closed form for the Stein discrepancy, which can be used as test statistic
in Monte Carlo procedures. 

As an illustration, \cite{fatima2025pure} assess the fit of a stochastic blockmodel to networks from the collection of Lazega's lawyers' networks \citep{lazega2001collegial}. These networks are derived from  a study on relationships between 71 attorneys (partners and associates)
in the three offices of a Northeastern US corporate law firm, during 1988--1991. The data provide  information on  advice  and friendship relationships among the 71 attorneys; they also contain information on attributes, such as the type of practice. From these data, several undirected networks are constructed, with lawyers as vertices. In the  {\emph{advice network}}, an edge between two lawyers indicates that they consult each other in a professional capacity, and an edge in the {\emph{friendship network}} indicates that the two lawyers socialise with each other outside work. For the advice network, at the 5\% level, the Monte Carlo test does not reject the null hypothesis of a Bernoulli random graph. For the friendship network, when assigning block membership according to the type of practice or to gender, a stochastic blockmodel is rejected. However, when assigning blocks according to which of the offices the lawyer works in, a stochastic blockmodel is not rejected.

In \cite{fatima2025pure}, Stein's method is also used to obtain explicit upper bounds on the Wasserstein distance
between the suitably normalised test statistic and the standard normal distribution.

\section{Exponential random graphs} \label{section_ergm}
In the previous section, we have concentrated on approximating the distributions of real valued statistics,
taking $\cX = \re$ in~\eqref{stein_identity}.  Here, we take~$\cX$ to be the space of all adjacency matrices, and approximate
the distribution of one random graph by another (simpler) distribution.  Suppose that $X$ and~$Y$ are random
elements of~$\bbA_n$, the set of
adjacency matrices on~$n$ vertices, and that their distributions $\law(X)$ and~$\law(Y)$ are each stationary distributions
of ergodic continuous time Markov processes $Z_X$ and~$Z_Y$ with generators $\cA_X$ and~$\cA_Y$, respectively.
If the Stein equation~\eqref{stein_equation} with
generator~$\cA_Y$ can be solved for $f = f_h$, for all~$h$ in a suitably large class of functions~$\cH$, and
if $\ex[\cA_Y f_h(X)]$ is uniformly small for all $h \in \cH$, then $\law(X)$ is correspondingly close to~$\law(Y)$,
Note that, typically, $f_h$ is given by 
\begin{equation}\label{recurrent-potential}
    f_h(a) \Eq - \int_0^\infty \ex\{h(Z_Y(t)) - \ex[h(Y)] \,|\, Z_Y(0) = a\}\,dt,\qquad a \in \bbA_n.
\end{equation}
Since also, from~\eqref{stein_equation}, $\ex[\cA_X f_h(X)] = \ex[h(X)] - \ex[h(X)] = 0$ for all
such~$h$, it follows that
\begin{equation}\label{generator_comparison}
     \ex[h(X)] - \ex[h(Y)] \Eq \ex[\cA_Y f_h(X)] - \ex[\cA_X f_h(X)] \Eq \ex[(\cA_Y - \cA_X) f_h(X)].
\end{equation}
Thus, the $d_\cH$-distance between $\law(X)$ and~$\law(Y)$ can be bounded in terms of the difference between
the two generators $\cA_Y$ and~$\cA_X$.

As an example, based on \cite{reinert2019approximating}, we consider approximating the distribution of an
exponential random graph model (ERGM) by a (much simpler) Bernoulli random graph.
\ignore{
Here, the approximation is not by a Poisson or a normal approximation, but instead by the distribution of a Bernoulli random graph. First, we consider an approach called comparison of Stein operators. We aim to estimate between the distance between two probability distributions $X$ and $Y$ where we have access to both Stein operators $\mathcal{A}_X$ and $\mathcal{A}_Y$. Due to the Stein identity \eqref{stein_identity} we can write
\begin{align*}
    \vert \mathbb{E}[h(X)-h(Y)] \vert = \vert \mathbb{E}[ \mathcal{A}_Y f_h(X) ] \vert = \vert \mathbb{E}[ ( \mathcal{A}_Y - \mathcal{A}_X ) f_h(X) ] \vert.
\end{align*}
The remainder of this section is based on \cite{reinert2019approximating} where this approach was used to construct bounds
on the difference of expectations between the Bernoulli and the ERGM. For two graphs $G=(V,A)$, $G'=(V',A')$ we let $t(A,A')$ the number edge-preserving injections from $V$ to $V'$. An injection $\sigma: V \rightarrow V'$ is edge-preserving if  $a_{u,w}=1$ implies $a_{\sigma(u),\sigma(v)}=1$.
}
An ERGM on~$n$ vertices, based on the counts of small connected subgraphs $G_1,\ldots,G_k$ in $(V,A)$, where $G_i$ has~$m_i$ vertices
and~$e_i$ edges, has probability distribution given quasi-explicitly by
\begin{align} \label{definition_ergm}
    \mathbb{P}(A=a) \Eq \cZ(\beta)^{-1}\exp\bigg( \sum_{i=1}^k \beta_i t_i(a) \bigg),\qquad a \in \bbA_n.
\end{align}
Here, $\beta = (\beta_1, \ldots, \beta_k) \in \mathbb{R}^k$ is a vector of parameters,
\[
   t_i(a) \Def t_{G_i}(a) \Def \frac{T_{G_i}(a)}{n(n-1)\cdots(n - m_i + 1)} \quad \mbox{and}\quad T_{G_i}(a) \Def
      \text{\# copies of } G_i \text{ in }a,
\]
and~$\cZ(\beta)$ is a (typically intractable) normalizing constant.
\ignore{
and
$t_i(a)$ denotes the number of copies ofand for graphs $G_1 = (V_1,A_1), \ldots, G_k = (V_k,A_k) $ we let
\begin{align*}
    t_i(a) = \frac{t(A_i,a)}{n(n-1) \ldots (n- \vert V_i \vert + 3)}
\end{align*}
and $n$ is the number of vertices with respect to the adjacency matrix $a$. Note that $t_i$ is defined on the set of all adjacency matrices with $n$ vertices.
}
The graph~$G_1$ is taken to be a single edge. \par

A suitable stochastic process~$Z_Y$ for use with the generator approach
is the \textit{Glauber dynamics} Markov chain. At the points of a Poisson process of rate~$1$,
a pair of vertices is chosen uniformly from~$V$, and the edge between the vertex pair is resampled
to be present or absent, conditional on the configuration on the rest of the graph.  This conditional
distribution can be written down from~\eqref{definition_ergm}, without needing to know~$Z(\beta)$. If the pair $\{u,v\}$
is chosen, and if~$p_1^{u,v}(a)$ denotes the conditional probability that the edge $\{u,v\}$ is present in $(V,A)$,
with~$p_0^{u,v}(a)$ the complementary probability, then 
\[
  p_1^{u,v}(a) \Def \frac{\exp\bigl( \sum_{i=1}^k \beta_i t_i(a_1^{u,v}) \bigr) }
   {\exp\bigl( \sum_{i=1}^k \beta_i t_i(a_1^{u,v}) \bigr)  + \exp\bigl( \sum_{i=1}^k \beta_i t_i(a_0^{u,v}) \bigr)};
   \qquad p_0^{u,v}(a) = 1 - p_1^{u,v}(a),
\]
for $a \in \bbA_n$, where $a_0^{u,v}$ resp.~$a_1^{u,v}$ denotes the adjacency matrix~$a$ with the value at $(u,v)$
(and therefore also at $(v,u)$) replaced by $0$ resp.~$1$.
\ignore{conditional on the present or absent edge of the chosen vertex pair according to \eqref{definition_ergm}.}
The infinitesimal generator of~$Z_Y$ is thus given by
\begin{equation} \label{ergm_stein_operator}
    \cA_Y f(a) \Eq \frac{1}{\binom n2} \sum_{\{u,v\}} \big( p_1^{u,v}(a) (f(a_1^{u,v}) - f(a))
        + p_0^{u,v}(a) (f(a_0^{u,v}) - f(a))  \big),
\end{equation}
where we sum over all vertex pairs $\{u,v\}$.
The differences $f_h(a_1^{u,v}) - f_h(a_0^{u,v})$ can be bounded using~\eqref{recurrent-potential}, together with known facts
about the Glauber dynamics starting from the states $a_1^{u,v}$ and~$a_0^{u,v}$.
\ignore{
$ \Delta_{(u,v)} f(a) = f(a_1^{(u,v)}) - f(a_0^{(u,v)})$
and $a_0^{(u,v)}$ resp.\ $a_1^{(u,v)}$ is the adjacency matrix $a$ with the value at $(u,v)$ (and therefore also at $(v,u)$)
equal to $0$ resp.\ $1$. This operator is a Stein operator for the ERGM. A solution to the Stein equation \eqref{stein_equation} for any real-valued function $h$ is 
given by
\begin{align*}
    f_h(a) = - \int_0^{\infty} \mathbb{E}\big[ h(A(t)) - \mathbb{E}[h(A) \, \vert \, A(0) = a] \big] dt, 
\end{align*}
where $A$ follows the ERGM and $A(\cdot)$ is the corresponding Glauber dynamics Markov chain.
}
An analogous Stein operator for $\cG(n,p)$ is given by
\begin{align*} 
    \cA_X f(a) \Eq \frac{1}{\binom n2} \sum_{\{u,v\}} \big( p(f(a_1^{u,v}) - f(a))  + (1-p)(f(a_0^{u,v}) - f(a))\big),
    \quad a \in \bbA_n.
\end{align*}

\ignore{In \cite{reinert2019approximating}, bounds on the differences $\Delta_s f(x) = f(x^{(s,1)}) - f(x^{(s,0)})$ are obtained by comparing the Glauber dynamics with two different starting points, $x^{(s,1)}$ and $x^{(s,0)}$. Here $x^{(s,1)}$ denotes the collection of edge indicators with $x_s$ set to equal 1, and $x^{(s,0)}$ is defined analously.}
The following theorem can then be deduced from~\eqref{generator_comparison}.  To state it, we define
\begin{align*}
    \Phi(p) \Def \sum_{i=1}^k \beta_i e_i p^{e_i-1}; \quad |\Phi| (p) \Def \sum_{\ell=1}^k |\beta_\ell| e_\ell p^{e_\ell-1}
    \ \text{ and }  \ \phi(p) \Def \frac{1+\tanh(\Phi(p))}{2}\,,
\end{align*}
and assume that $p^{\star} \in [0,1]$ satisfies $p^{\star} = \phi(p^{\star})$.
For any $g\colon \bbA_n \to \re$, we 
write 
$
  \Delta_{\{u,v\}}g(a) := g(a_1^{u,v}) - g(a_0^{u,v}), 
$  
noting that, for $X \sim \cG(n,p)$, the distribution of $\Delta_{\{u,v\}}g(X)$ is the same for all~$\{u,v\}$.
Finally, we set $\|\Delta g\| := \sup_{\{u,v\} \subset V } \sup_{a \in \bbA_n} \vert \Delta_{(u,v)}g(a) \vert$,
and define the class~$\cH_n$ of Lipschitz functions $\cH_n := \{h\colon \bbA_n \to \re; \|\Delta h\| \le 1/\binom n2\}$.
Such functions include the normalised subgraph counts~$t_H$, which
feature prominently in notions of dense graph convergence \citep{borgs2008convergent, borgs2012convergent}.
\begin{theorem} \label{theorem_exponential_random graph}
Suppose that~$Y$ is sampled from an ERGM with probabilities as in~\eqref{definition_ergm}, and that $X \sim \cG(n,p^{\star})$.
Then, with the above definitions, in the subcritical case $|\Phi|'(1) < 2$, it follows that
\ignore{
Assume that $p^{\star} \in [0,1]$ satisfies $p^{\star} = \phi(p^{\star})$ and let $P^{\star}= \max\{p^{\star}, 1-p^{\star} \}$. Moreover, we set
    \begin{gather*}
        \alpha_1 = \frac{1}{2} \big( \Phi'(p^{\star}) + P^{\star} \Phi''(1) \big), \\
        \alpha_2 = \phi'(p^{\star}) + \frac{1}{2} \bigg( \frac{4}{3\sqrt{3}} (P^{\star} + n^{-1}) \Phi'(1) \big( \Phi'(p^{\star}) + P^{\star} \Phi'(1) \big) + P^{\star} \Phi''(1) \mathrm{sech}^2(\Phi(p^{\star})) \bigg)
    \end{gather*}
}%
    \begin{align*}
        d_{\cH_n}(\law(Y),\law(X))
        \Le   \frac1{2(2 - |\Phi|'(1))}\,
        \sum_{i=2}^k |\beta_i |
            \sqrt{\mathbb{V}[\Delta_{(u,v)} t_i(X)]}.
    \end{align*}
\end{theorem}
\nin Using arguments from \cite{rucinski1988small}, it can be shown that, for any fixed~$p$, the random variables
$\Delta_{(u,v)} t_i (X)$ have variance of order $1/n$ when $X\sim \mathcal{G}(n, p)$. Hence, if $|\Phi|'(1) < 2$,
it follows that $d_{\cH_n}(\law(Y),\law(X))$ is of order~$O(n^{-1/2})$.
\ignore{
A particular group of functions  of interest are the scaled
subgraph counts~$t_H(a)$ of small connected graphs~$H$.  These
 feature prominently in notions of graph convergence (see for example \cite{borgs2008convergent, borgs2012convergent}). For such functions~$h$,
 the bound in Theorem~\ref{theorem_exponential_random graph} can be seen to be of order $O(n^{-1/2})$.}

We give two examples of an application of Theorem \ref{theorem_exponential_random graph} for $k=2$.
If~$G_2$ is a two star and  $\vert \beta_2 \vert <1$, then there is a unique $p^{\star}$ satisfying $p^{\star} = \phi(p^{\star})$, and
\begin{align*}
    d_{\cH_n}(\law(Y),\law(X))
       \Le  \frac{\vert \beta_2 \vert}{4(1-\vert \beta_2 \vert)}\, \frac{\sqrt{8p^{\star}(1-p^{\star}) }  }{\sqrt{n-2}}\,.
\end{align*}
If~$G_2$ is a triangle and  $\vert \beta_2 \vert < 1/3$, then there is a unique $p^{\star}$ satisfying $p^{\star} = \phi(p^{\star})$, and
\begin{align*}
    d_{\cH_n}(\law(Y),\law(X))
        \Le \frac{\vert \beta_2 \vert}{4(1-3 \vert \beta_2 \vert)} \frac{6p^{\star}\sqrt{1-(p^{\star})^2 }  }{\sqrt{n-2}}\,.
\end{align*} 

In \cite{xu2021stein}, the authors propose a goodness-of-fit test for the ERGM, based on the Stein discrepancy~\eqref{Stein_discrepancy}
derived from the operator \eqref{ergm_stein_operator}. Theorem \ref{theorem_exponential_random graph}
is used to assess the error incurred when test statistics evaluated at an ERGM are treated as having come from the corresponding
Bernoulli model $\cG(n,p^\star)$.

We conclude the section with recent results regarding the approximation of count statistics in ERGMs. In \cite{fang2025normal} the authors employ Stein's method by generalising an approach from \cite{chatterjee2008new} to obtain a central limit theorem for functionals of nonlinear exponential families of the form
\begin{align*}
    q(y) = \frac{\exp ( h(y) ) g(y)}{\mathbb{E}[\exp ( h(X) )]}, \qquad y \in \mathbb{R}^d, 
\end{align*}
where $X \in \mathbb{R}^d$ is a random vector of independent random variables with density function $g$. The result is then applied to the ERGM to obtain normal approximations of suitably normalised subgraph counts (including the number of edges as a special case). \cite{fang2025conditional} then extended these results to subgraph counts conditional on the number of edges and \cite{winstein2025quantitative} to settings beyond the subcritical regime.

\section{Spatial random graphs} \label{section_spatial_graphs}

A geometric graph consists of a point configuration~$\xi$, together with a set of edges between the points.
In a simple model of a random geometric graph, $\xi$ consists of~$n$ points, uniformly distributed in $[0,1]^2$, and,
for some $r > 0$, any
two points less than distance~$r$ apart are joined by an edge; see for example \cite{penrose2003}.
Much more complicated random models can be envisaged,
with~$\xi$ sampled from a general point process, and with edges assigned according to an arbitrary joint distribution
depending on~$\xi$.  The purpose of this section is to show how to use Stein's method to quantify the error involved in
approximating a complicated model by a simpler one.  We mainly base our discussion on
\cite{schuhmacher2024stein} and the PhD thesis \cite{wirth2025spatial}.
\ignore{
where bounds for the distance to a generalised random geometric graph are obtained with respect to integral probability metrics. A generalised geometric graph can be represented by a counting measure $\xi$ on the underlying space. Here again a key task is to approximate the model by a simpler model. This now requires a framework for comparing two counting measures on the same space.
}

We assume, for convenience, that a configuration~$\xi$ consists of {\it distinct\/} points belonging to a compact set
$\bbX \subset \mathbb{R}^d$. Writing $\mathfrak{R} = \mathfrak{R}(\bbX)$ for the space of all finite counting measures on~$\bbX$, we note that a configuration~$\xi$ can be interpreted as such a measure, with~$\xi(B)$ the number of points of~$\xi$ lying
in $B \subset \bbX$.
Given a configuration~$\xi$, the set of edges~$\sigma$ between the pairs of points of~$\xi$ is chosen according to
a probability distribution~$Q(\xi,\cdot)$ (the \textit{edge kernel}) on~$\xi^{\langle 2 \rangle}$, the set two point 
subsets of~$\xi$.  For convenience, we represent~$\xi^{\langle 2 \rangle}$ as an element of
$\mathfrak{R}_2 = \mathfrak{R}(\bbX^{\langle 2 \rangle})$, the set of counting measures on the space
$\bbX^{\langle 2 \rangle}$ of two element subsets of~$\bbX$, and
take~$Q(\xi,\cdot)$ to be defined on the whole of~$\mathfrak{R}_2$, albeit concentrated on
$\mathfrak{R}(\xi^{\langle 2 \rangle})$, to which~$\sigma$ belongs. The pair $(\xi,\sigma)$ constitutes a realisation of a
random graph, whose probability measure on $\mathfrak{R} \times \mathfrak{R}_2$ is obtained
\ignore{
Moreover, we let $ X^{\langle 2 \rangle} = \{ (x,y) \in \bbX^2 \, \vert \, x < y \}$ denote the set of all two-element subsets of $\bbX$ and accordingly $\mathfrak{R}_2 = \mathfrak{R}(\bbX^{\langle 2 \rangle})$.  We call a probability kernel $Q$ an \textit{edge kernel} if the measure $Q(\xi, \cdot)$ for $\xi \in \mathfrak{R}$ is concentrated on the set $\mathfrak{R}_2 \vert_{\xi^{\langle 2 \rangle}}$.
}
from a probability distribution  $P$ on~$\mathfrak{R}$ and an edge kernel $Q$ on~$\mathfrak{R}_2$, and is denoted by
$P \otimes Q$:
the  expectation of any bounded measurable function
$f\colon \mathfrak{R} \times \mathfrak{R}_2 \rightarrow \mathbb{R}$ with respect to $P \otimes Q$
is defined to be
\begin{align*}
    P \otimes Q(f) \Eq \int_{\mathfrak{R}} \int_{\mathfrak{R}_2} f(\xi, \sigma) Q(\xi, d\sigma) P(d\xi).
\end{align*}
We write $(\Xi,\Sigma) \sim P \otimes Q$
for a random element of $\mathfrak{R} \times \mathfrak{R}_2$ with this distribution, and
\begin{align}\label{graph_space}
    \mathbb{G} =\{ (\xi, \sigma) \, \vert \, \xi \in \mathfrak{R}, 
           \sigma \in \mathfrak{R}(\xi^{\langle 2 \rangle})\} \}
\end{align}
for the space of all possible realisations.
Note that, to obtain configurations of distinct points, we require~$P$ to assign probability zero to the set of all
counting measures in~$\mathfrak{R}$ that have atoms greater than~$1$.
\ignore{
we write $(\Xi,\Sigma) \sim P \otimes Q$. Furthermore, we write 
$\mathbb{G} =\{ (\xi, \sigma) \, \vert \, \xi \in \mathfrak{R}, 
\sigma \in \mathfrak{R}_2 \vert_{\xi^{\langle 2 \rangle}} \}$ where we wrote 
$\mathfrak{R}_2 \vert_{\xi^{\langle 2 \rangle}}$ for the space
$\mathfrak{R}_2 $ restricted with respect to the points in $\xi$. 
}

\ignore{
We introduce the random geometric graph. For this purpose we recall that a point process is a (possibly infinite) set of random points $X_1, X_2, \ldots \in \mathbb{R}^d$ which can be identified as a random measure $\Xi = \sum_i \delta_{X_i}$  such that $\Xi(B) < \infty$ almost surely for any compact set $B \in \mathbb{R}^d$.
}
A Poisson point process~$\Xi \sim \Pop(\bla)$ on~$\bbX$, with intensity measure $\bla$ on~$\bbX$,
is a random element of~$\mathfrak{R}$ such that
\begin{itemize}
    \item the number of points $\Xi(B)$ in any measurable set $B \subset \bbX$ satisfies $\law(\Xi(B)) = \Po(\bla(B))$, and
    \item for disjoint measurable sets $B_1, \ldots, B_n \subset \bbX$, the numbers $\Xi(B_1), \ldots, \Xi(B_n)$ are
        independent.
\end{itemize}
The independence edge kernel, based on the connection probabilities
$\kappa\colon \bbX^{\langle 2 \rangle} \rightarrow [0,1]$, is defined by
\begin{align*}
    Q^{\kappa} (\xi, \{ \sigma \}) \Eq \prod_{\{x,y \} \in \xi^{\langle 2 \rangle}}
                 \kappa(x,y)^{\sigma(\{x,y\} )} (1- \kappa(x,y))^{1-\sigma(\{x,y\} ) },\qquad 
                 \sigma \in \mathfrak{R}(\xi^{\langle 2 \rangle});
\end{align*}
conditional on the positions~$\xi$ of the points,
edges are assigned between pairs of them independently, with probabilities depending on their positions.
A {\it (soft) random geometric graph\/} $\mathrm{RGG}(\bla, \kappa)$ 
with (atomless) intensity measure~$\bla$ and connection probabilities~$\kappa$ is defined to have the distribution 
$\Pop(\bla) \otimes Q^{\kappa}$. 
\ignore{
For the construction of the actual graph we interpret $\xi$ as the vertex and $\sigma$ as the edge process, i.e.\ a realisation $(\xi, \sigma)$ gives access to the position of the vertices in $\mathcal{X}$ and the edges between those vertices whereby an edge between two vertices $x$ and $y$ occurs with probability $\kappa(x,y)$. 
}
A typical example of the function~$\kappa$ is
\begin{align*} 
    \kappa(x,y) \Eq \begin{cases} 1, & \Vert x-y \Vert \leq r; \\  0, & \Vert x-y \Vert > r,  \end{cases}
\end{align*}
where $r > 0$ is the threshold; two points $x$ and $y$ are joined by an edge if they are at distance at most~$r$ from
one another. Another example is given by the \textit{Rayleigh fading activation function} (see for example \cite{duchemin2023random}), defined by
\begin{align} \label{kappa_threshold-fading}
    \kappa(x,y) \Eq \exp\bigg( - \zeta \bigg( \frac{\Vert x - y \Vert}{r} \bigg)^{\eta} \bigg),
\end{align}
where $r, \eta, \zeta > 0$ are parameters.

\begin{figure}[ht]
    \centering
    \begin{minipage}{0.48\textwidth}
        \centering
        \includegraphics[width=\linewidth]{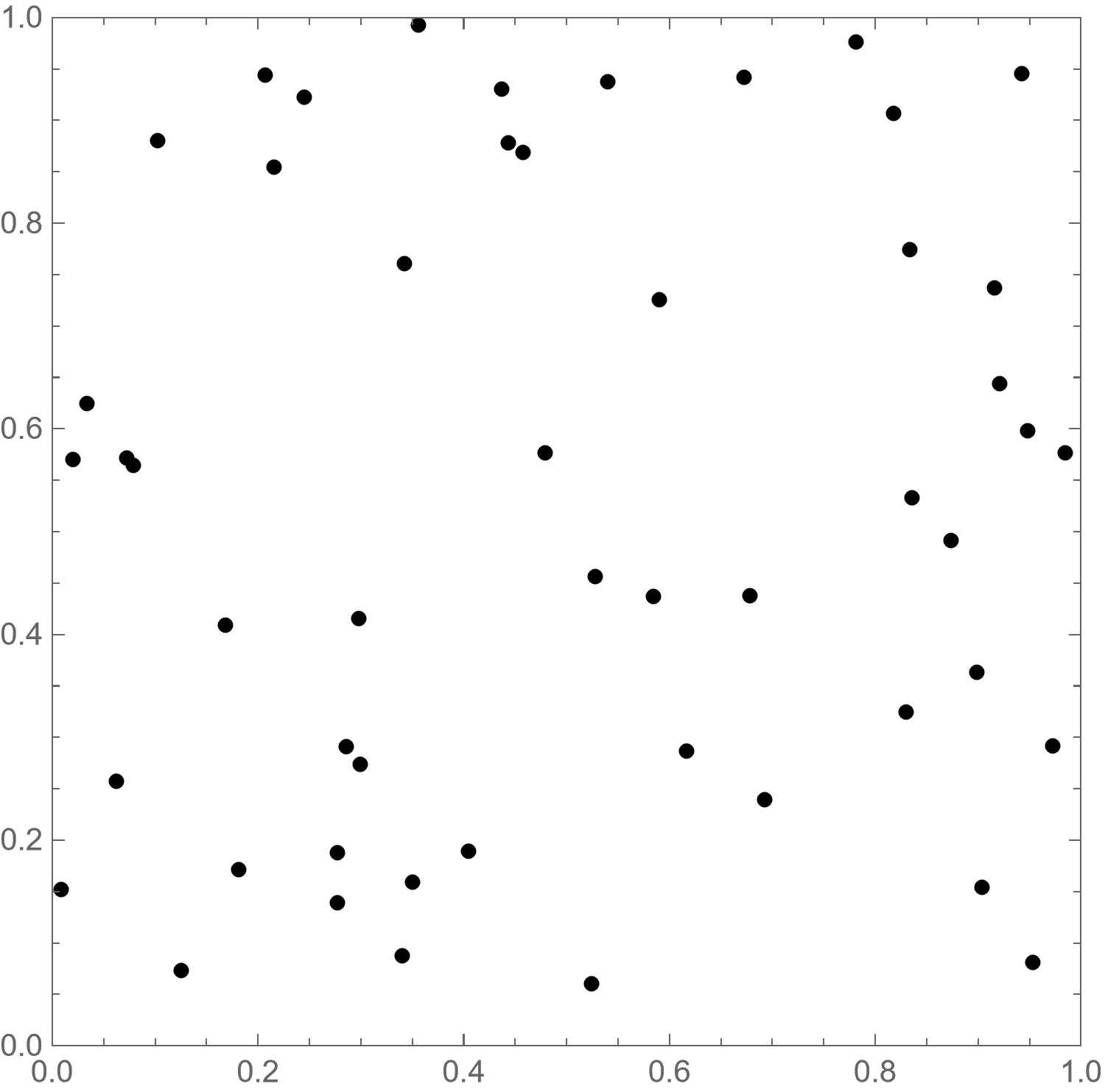}
        \caption{Realisation of a Poisson point process with $\bla$ the Lebesgue measure on $\mathbb{R}^2$ 
        with intensity $50$ and $\bbX=[0,1]^2$.}
    \end{minipage}
    \hfill
    \begin{minipage}{0.48\textwidth}
        \centering
        \includegraphics[width=\linewidth]{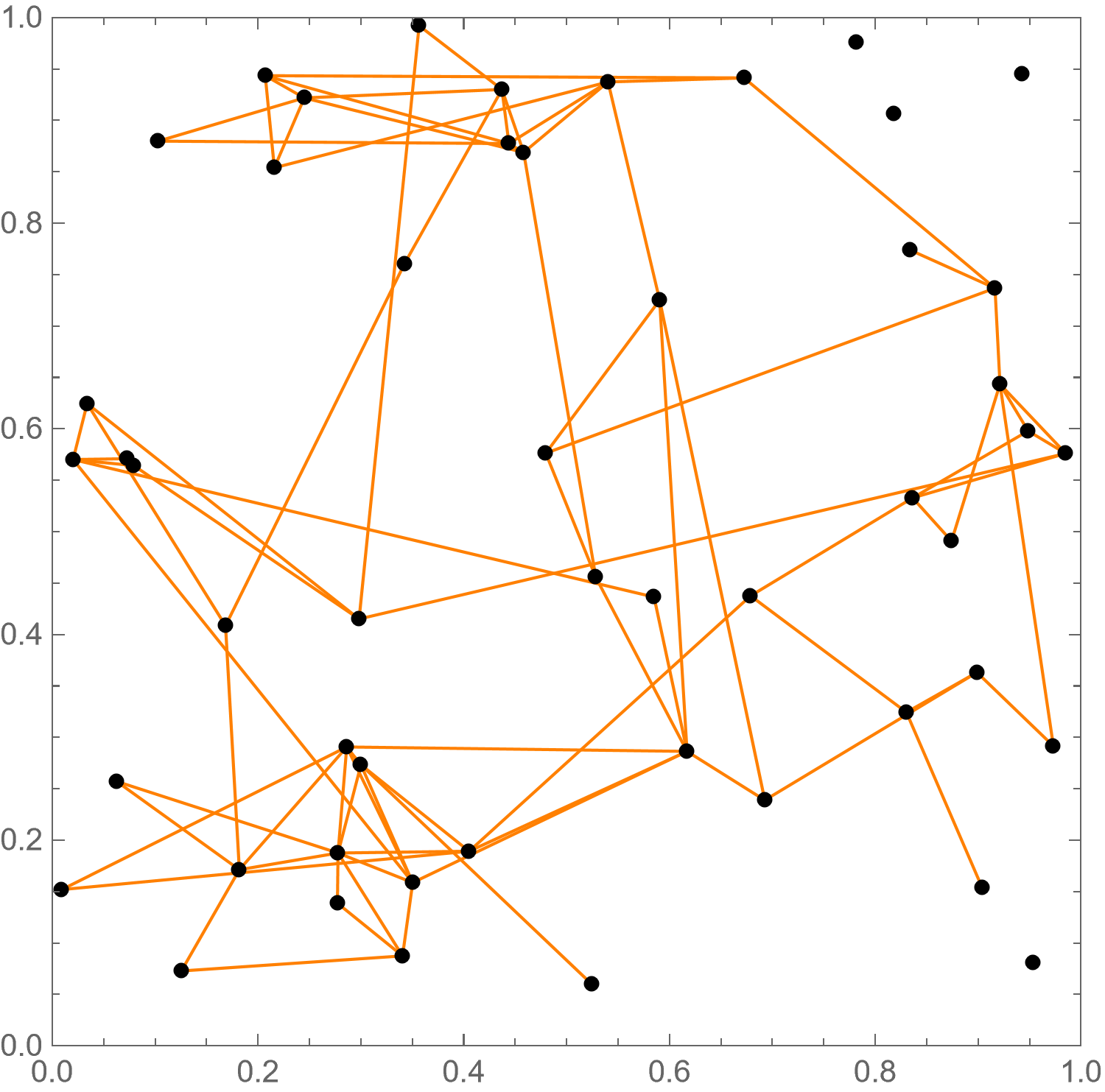}
        \caption{Realisation of a random geometric graph on the same set of points, with~$\kappa$ as in
        \eqref{kappa_threshold-fading} with $\zeta=8$ and $\eta=r=1$.}
    \end{minipage}
\end{figure}

A Stein operator for $\mathrm{RGG}(\bla, \kappa)$ can then be constructed using the generator approach. As remarked in Section~\ref{Poisson_approximation}, the Stein operator~$\tcA_\la$ 
for the Poisson distribution~$\Po(\la)$, given in~\eqref{Poisson_Stein_operator2},
can be recognized as the 
generator of an immigration-death process with immigration rate~$\lambda $ and unit per capita death rate.
\ignore{
First, let us consider the following birth-and-death process $Z_t, t \geq 0$ which takes values in $\mathbb{N}$: Let $Z(t)$ be the number of individuals present at time $t$ where individuals arrive according to a (one-dimensional) Poisson process with constant rate $\lambda>0$. Once they arrived, they have lifetimes independently of each other and the birth process which are exponentially distributed with mean $1$. Therefore, given that $Z(t)=j$, we wait for an $Exp(\lambda + j)$ amount of time until the next event which is a birth with probability $\lambda/(j+\lambda)$ and a death with probability $j/(j+\lambda)$. This process is then Markov with stationary distribution $P(\lambda)$. The latter
}
This representation can be generalised to the random geometric graph, in the form of a graph-valued birth-and-death 
process~$Z_{\bla,\kappa}$ with state space~$\mathbb{G}$, as defined in~\eqref{graph_space}. 
Given $Z_{\bla,\kappa}(t) = (\xi,\sigma)$, the transition rate for adding a point is~$\bla(\bbX)$,
and for deleting a point is~$|\xi|$.  If a point is added, it is assigned its position according to the measure 
$\bla (\cdot) /\bla(\bbX)$, and edges are added between the new point~$x$ and the old points~$\xi$ according to
\begin{align}\label{Q2-def}
    Q_2^{\kappa} (\xi, x; \{\sigma_2 \} ) 
      = \prod_{y \in \xi} \kappa(x,y)^{\sigma_2(\{x,y\})} (1- \kappa(x,y))^{1-\sigma_2(\{x,y\} ) },\quad
          \sigma_2 \in \mathfrak{R}(\{x,y\}\colon y \in \xi),
\end{align}
thus independently of the existing edge structure $\sigma$, but not of the positions of the points~$\xi$. 
If a point is deleted, it is chosen uniformly among all points of~$\xi$, and all edges incident upon it
are also deleted.  As in \cite[Section~4.2 and Section~4.11, Problem 5]{ethier2009markov}, 
it can be shown that the generator of~$Z_{\bla,\kappa}$, acting on $f\colon \bbG \rightarrow \mathbb{R}$, is given by
\begin{align*}
    \cA^{\bla,\kappa} f(\xi, \sigma) \Eq 
       \int_{\bbX} \mathbb{E}[f(\xi + & \delta_x , \sigma + T_{\xi, x}) - f(\xi, \sigma)] \bla(dx)  \\ 
    & + \int_{\bbX} f(\xi - \delta_x , \sigma_{\xi,x} ) - f(\xi, \sigma)
       \xi(dx) ,
\end{align*}
where $ T_{\xi, x} \sim  Q_2^{\kappa} (\xi, x; \cdot )$ and $\sigma_{\xi,x}$ denotes the restriction of~$\sigma$
to $(\xi-\delta_x)^{\langle 2 \rangle}$. Then, given a bounded measurable function 
$h\colon \mathbb{G} \rightarrow \mathbb{R}$, the Stein equation
\begin{align*} 
    h(\xi, \sigma) - \mathbb{E}[h(H,T)] \Eq \cA^{\bla,\kappa} f(\xi, \sigma),
\end{align*}
where $(H,T) \sim \Pop(\bla) \otimes Q^{\kappa}$, can be solved by 
\begin{align}\label{PP_stein_solution}
    f_h(\xi, \sigma) \Eq - \int_0^{\infty} \mathbb{E}[f((H_s,T_s)^{(\xi, \sigma)} ] - \mathbb{E}[f(H,T)]\, ds,
\end{align}
where $((H_s,T_s)^{(\xi, \sigma)},\,s \geq 0)$ denotes the Markov process~$Z_{\bla,\kappa}$, started in 
$(\xi, \sigma)$. 

In \cite{schuhmacher2024stein}, the main interest lies in approximating the distribution of  
random elements of~$\bbG$ whose distributions admit a so-called \textit{Papangelou kernel} (see for example 
\cite[Chapter~8]{kallenberg2017random}).  Their approximations are carried out in terms of 
\textit{generalised random geometric graphs}, in which the point pattern is a realisation of a point process 
that admits a density 
with respect to a Poisson point process~$\Pop(\bla)$ (a so-called \textit{Gibbs process}). For such processes,
the generator~$\cA^{\bla,\kappa}$ is replaced by something rather more complicated:  for instance, the rate of 
creation of new points depends not only on the
position of the new point, but also on the existing configuration~$\xi$.

Here, we give just a flavour of their general results, 
presenting a simplified version of their main theorem; we compare the distributions of two random geometric 
graphs with different 
vertex and edge processes. 
To state it, we define
\begin{align*}
    \Delta_V^{\kappa} f(\xi, \sigma) &:= \sup_{x \in \bbX \setminus \xi}  
       \big\vert \mathbb{E}[f(\xi+\delta_x,\sigma+T_{\xi,x})- f(\xi,\sigma) ] \big\vert, 
       \quad \mbox{for}\ T_{\xi,x} \sim Q_2^\kappa;\\
    \Delta_E f(\xi, \sigma) &:= 
       \sup_{\substack{x \in \bbX \setminus \xi,\, y \in \xi \\ 
          \sigma_2 \in \mathfrak{R}(\{x,z\}\colon z \in \xi \setminus \{y\}) }}
    \bigl\vert f(\xi+\delta_x,\sigma+\sigma_2+\delta_{\{x,y\}})- f(\xi+\delta_x,\sigma + \sigma_2 ) \bigr\vert.
\end{align*}
The quantity
$\Delta_V^{\kappa} f(\xi, \sigma)$ can be interpreted as the maximal expected change in~$f$,
with $T_{\xi,x} \sim Q_2^\kappa$ as in~\eqref{Q2-def},
when a point (and its edges) are added to the graph $(\xi,\sigma)$;
$\Delta_E f(\xi, \sigma)$ represents the maximal difference in~$f$ if, when a point~$x$ and its edges are added to the
graph $(\xi,\sigma)$, a particular edge between $x$ and a point of~$\xi$ is present or absent.

\begin{theorem} \label{theorem_ipm_geom_random_graphs}
    Let $(H,T) \sim \Pop(\bla_1) \otimes Q^{\kappa_1}$ and $(\Xi,\Sigma) \sim \Pop(\bla_2) \otimes Q^{\kappa_2}$ 
    be random geometric graphs. 
    Define
\[
     \Delta\kappa(\xi,x) \Def \sum_{y \in \xi}|\kappa_1(x,y) - \kappa_2(x,y)|.
\]
    Then, for all $h\colon \mathbb{G} \rightarrow \mathbb{R}$ bounded and measurable,
    \begin{align*}
        \vert \mathbb{E}[h(&\Xi,\Sigma)] - \mathbb{E}[h(H,T)] \vert \\
        &\Le  \mathbb{E} \biggl[  \int_{\bbX} \Delta_V^{\kappa_2} f_h(\Xi, \Sigma) \vert \bla_1(dx) 
                  - \bla_2(dx) \vert \biggr] 
         + \mathbb{E} \biggl[  \int_{\bbX} \Delta_E f_h(\Xi, \Sigma) \Delta\kappa(\Xi,x)  \bla_1(dx)  \biggr].
    \end{align*}
\end{theorem}

\ignore{
In the sequel we define the remaining expressions and briefly give some intuitions behind both terms from the bound in Theorem \ref{theorem_ipm_geom_random_graphs}. 
The first term measures the difference of the vertex processes $\Xi$ and $H$, i.e. the difference in terms of the intensity measures $\lambda$ and  $\tilde{\lambda}$. Regarding the second term, we define
\begin{align*}
    \Vert \hat{v}(\xi,x)) \Vert_1 = \sum_{i=1}^{\vert \xi \vert} \Big\vert \kappa(x_i,x) - \tilde{\kappa}(x_i,x) \big) \Big\vert,
\end{align*}
where $\xi = (x_1,x_2,\ldots)$, i.e. $\Vert \hat{v}(\xi,x) \Vert_1$ compares the probabilities for the presence of an edge between an additional vertex $x$ and the vertices of an existing graph $(\xi, \sigma)$ regarding the edge processes $\Sigma$ and $T$. 
}

For a general bounded measurable function~$h\colon \bbG \to \re$, it is easy to deduce from~\eqref{PP_stein_solution} that
$\Delta_V^{\kappa} f(\xi, \sigma) \le 2\|h\|_\infty$ and that $\Delta_E f(\xi, \sigma) \le \|h\|_\infty$.
Theorem~\ref{theorem_ipm_geom_random_graphs} then implies that
\begin{align*}
    \vert \mathbb{E}[h(&\Xi,\Sigma)] - \mathbb{E}[h(H,T)] \vert \\ &\Le 
       \|h\|_\infty \bigg( 2\int_{\bbX} |\bla_1(dx) - \bla_2(dx)| 
          + \int_{\bbX} \int_{\bbX} |\kappa_1(x,y) - \kappa_2(x,y)|\bla_1(dx)\bla_2(dy) \bigg),
\end{align*}
in which the impact of the differences between the pairs $(\bla_1,\kappa_1)$ and $(\bla_2,\kappa_2)$ is clearly visible.
Taking all functions~$h$ having $\|h\|_\infty \le 1/2$ yields a corresponding bound on 
$d_{\TV}\bigl(\Pop(\bla_1) \otimes Q^{\kappa_1},\Pop(\bla_2) \otimes Q^{\kappa_2}\bigr)$.

For many applications, total variation distance is too strong, and Wasserstein metrics, based on 
classes of test functions~$h$ 
that are Lipschitz with respect to choices of the underlying topology on~$\bbG$, are more useful.  
The intuition is that changing the positions of points by only a small amount, or changing a small
number of edges, should not greatly affect the value of~$h$. For
such classes of functions, the quantities $\Delta_V^{\kappa} f(\xi, \sigma)$ and~$\Delta_E f(\xi, \sigma)$
may be a lot smaller, leading to more precise bounds: see \cite{schuhmacher2024stein} for much more
on the subject.

As for Bernoulli models, limit theorems for functions of random geometric graphs are available; see for example \cite{schulteyukich}. In applications, the spatial coordinates of the points in a random geometric graph may not be given, with the only information consisting of the adjacency matrix. In this setting, the issue of estimating the underlying latent space structure arises, as in \cite{hoff2002}; \cite{sosa} gives a review of this field. \par

This short survey provides some sketches of the use of Stein's method in network analysis; the literature cited contains much more information about what is possible.

\backmatter

\bmhead{Acknowledgements}
AF is funded by EPSRC Grant EP/T018445/1. GR is funded in part by EPSRC grants EP/T018445/1, EP/V056883/1, EP/Y028872/1,  and EP/X002195/1.

\section*{Declarations}

For the purpose of open access, the authors have applied a Creative Commons Attribution (CC BY) licence to any Author Accepted Manuscript version arising from this submission.



\end{document}